\documentclass[11pt]{amsart}
\usepackage{latexsym,amssymb,amsfonts,amsmath,graphicx,fullpage,url}
\usepackage[vcentermath,enableskew]{youngtab}
\usepackage{color}
\usepackage{cite}
\usepackage[all]{xy}
\usepackage{verbatim}
\newtheorem{theorem}{Theorem}[section]
\newtheorem{proposition}{Proposition}[section]
\newtheorem{lemma}[theorem]{Lemma}
\newtheorem{corollary}[theorem]{Corollary}
\theoremstyle{definition}
\newtheorem{definition}[theorem]{Definition}

\theoremstyle{remark}

\numberwithin{equation}{section}

\def\C{{\mathcal{C}}}
\def\S{{\mathcal{S}}_n}

\def\L{{\mathcal{L}}}
\def\H{{\mathcal{H}}}

\def\T{{\mathcal{T}}}

\def\R{\mathbb{R}}

 \def\K{K_P^{\rm root}}
\def\x{{\bf x}}

\def\G{\Psi_P(\x)}

\newcommand{\be}{\begin{equation}}
\newcommand{\ee}{\end{equation}}

\newcommand{\bd}{\begin{definition}}
\newcommand{\ed}{\end{definition}}
\newcommand{\bt}{\begin{theorem}}
\newcommand{\et}{\end{theorem}}
\newcommand{\bl}{\begin{lemma}}
\newcommand{\el}{\end{lemma}}
\newcommand{\bp}{\begin{proposition}}
\newcommand{\ep}{\end{proposition}}
\newcommand{\bc}{\begin{corollary}}
\newcommand{\ec}{\end{corollary}}

\newtheorem*{theorem1*}{Theorem \ref{thm:main}}

\newcommand{\old}[1]{}

\author[Karola\ M\'esz\'aros]{Karola M\'esz\'aros}
\email{karola@math.cornell.edu}
\address{
Department of Mathematics, Cornell University, Ithaca, NY 14853
}
\thanks{The author was partially supported by a National Science Foundation Grant  (DMS 1501059).}

\title{Calculating Greene's function via root polytopes and subdivision algebras}

\begin{document}

\begin{abstract}  Greene's rational function $\Psi_P(\x)$ is a  sum of certain rational functions in $\x=(x_1, \ldots, x_n)$ over the linear extensions of the poset $P$ (which has $n$ elements), which he introduced in his study of the Murnaghan-Nakayama formula for the characters of the symmetric group. In recent work Boussicault, F\'eray, Lascoux and Reiner showed that $\G$ equals a valuation on a cone and calculated $\G$ for several posets this way. In this paper we give an expression for   $\G$ for any poset $P$.  We obtain such a formula using dissections of root polytopes. Moreover, we use the subdivision algebra of root polytopes to show that in certain instances $\G$ can be expressed as  a product formula, thus giving a compact alternative proof of  Greene's original result and its generalizations.
\end{abstract}

 \maketitle

\section{Introduction}

Given a poset $P$ on the set $[n]=\{1, \ldots, n\}$, Greene's rational function is defined by 
\be \label{green} \G=\sum_{w \in \L(P)}w\left(\frac{1}{(x_1-x_2)(x_2-x_3)\cdots(x_{n-1}-x_n)}\right).
\ee

It was introduced by Greene \cite{greene} in his work on the Murnaghan-Nakayama formula. In \cite{reiner} Boussicault, F\'eray, Lascoux and Reiner showed that 
  
  \be \label{s} \G=s(K_P^{\rm root}; \x), \ee
  where \be K_P^{\rm root}=\R_+\{e_i-e_j \mid i<_P j\}=\R_+\{e_i-e_j \mid i\lessdot_P j\} \ee and \be s(K; \x):=\int_K e^{-\langle \x, v\rangle} dv,\ee for $K$ a polyhedral cone in a Euclidean space $V$ with inner product $\langle \cdot, \cdot \rangle$. 
  \medskip
 
Next we explain  two important results about calculating $\G$. Further work on $\G$ appeared in \cite{bf,b1, b2, i}.

\subsection{Greene's Theorem.} Let $P$ be a {\it strongly planar} poset, meaning that the Hasse diagram of $P \sqcup \{\hat{0}, \hat{1}\}$ has a planar embedding with all edges directed upward in the plane. For a strongly planar poset $P$ the edges of the Hasse diagram of $P$  dissect the plane into bounded regions $\rho$ such that the set of vertices of $P$ in the boundary of $\rho$ are two chains starting and ending at the same two elements, $min(\rho)$ and $max(\rho)$, respectively. Denote by $b(P)$ the set of bounded regions into which the  Hasse diagram of $P$  dissects the plane.

\vspace{.15in}

\noindent \textbf{Greene's Theorem.} \cite{greene} \textit{For any strongly planar poset $P$, }

\be \label{eq:greene} \G=\frac{\prod_{\rho \in b(P)}(x_{min(\rho)}-x_{max(\rho)})}{\prod_{i\lessdot_{P}j}(x_i-x_j)}.\ee

 \subsection{Boussicault's, F\'eray's, Lascoux's and Reiner's Theorem.} A beautiful theorem appearing in \cite{reiner} gives an expression for $\G$ in case of some posets $P$ whose Hasse diagrams are bipartite graphs in terms of certain lattice paths. The setup is as follows. Let $D$ be a skew Ferrers diagram in English notation, and let us labels its rows from top to bottom by $1, 2, \ldots, r$ and its columns from right to left by $1, 2, \ldots, c$.  With this labeling the northeasternmost point of $D$ is $(1,1)$ and the southwesternmost is $(r,c)$. The \textbf{bipartite poset} $P_D$ is a poset on the set $\{x_1, \ldots, x_r, y_1, \ldots, y_c\}$ with order relations $x_i<_P y_j$ if and only if $(i,j) \in D$.

  \vspace{.15in}

\noindent \textbf{BFLR Theorem.} \cite{reiner} \textit{For any skew diagram $D$, }

\be \label{eq:reiner}\Psi_{P_D}(\x)=\sum_{\pi}\frac{1}{\prod_{(i,j)\in \pi}(x_i-y_j)},\ee
\textit{where the sum runs over all lattice paths $\pi$ from $(1,1)$ to $(r,c)$ inside $D$ that take steps either one unit south or west.}

\medskip

\subsection*{Roadmap of the paper.}The objective of this paper is to (1) give a combinatorial expression of  $\G$ for any poset $P$, (2) give an alternative proof of the BFLR Theorem and (3) generalize Greene's theorem.  We accomplish (1) and (2)  in Section \ref{arb}, while we do (3) in  Sections \ref{sec:lift} and \ref{sec:gen}. In Sections \ref{sec:lift} and \ref{sec:gen} we also study the integer point transform of the root cone, which can be seen as a more refined invariant of the cone  than Greene's function. The   integer point transform of the root cone and   generalizations of Greene's theorem were also investigated in \cite{reiner}. Our tools will be root polytopes and their subdivision algebras, the latter of which was introduced in \cite{root1} and put to use in \cite{prod, mm, h-poly1, h-poly2, pipe1, pipe2}.

\section{Greene's function for an arbitrary poset}
\label{arb}
The purpose of this section is twofold. First we show how to express $\G$ for any poset $P$ in terms of $\G$ for posets $P$ whose Hasse diagrams are alternating graphs. Then we give an expression for $\G$ for a posets whose Hasse diagrams are alternating graphs, thereby also obtaining an expression for $\G$ for any poset $P$.  Finally,  we show that for certain posets $P$ whose Hasse diagrams are  bipartite graphs we can write $\G$ as a nice summation formula. The latter result originally appeared in the work of Boussicault, F\'eray, Lascoux and Reiner \cite{reiner} who used triangulations of order polytopes in their proof. We phrase our proof in terms of root polytopes. The point of view of this paper is that (dissections of) root polytopes (and the root cone) are the unifying approach to the calculation of $\G$.

A \textbf{root polytope} (of   type $A_{n-1}$) is  the convex hull of the origin and some of the points $e_i-e_j$ for $1\leq i<j \leq n$. Given a graph $G$ on the vertex set $[n]$ we associate to it the root polytope \be \tilde{Q}_G={\rm ConvHull}(0, e_i-e_j\mid (i,j) \in E(G), i<j).\ee 

It can be seen that $\tilde{Q}_G$ is a simplex if and only if $G$ is acyclic  and to emphasize this we sometimes denote $\tilde{Q}_G$ for acyclic graphs $G$ by $\tilde{\Delta}_G$.

The posets $P$ we work with in this section are on  the set $[n]$ and they are labeled naturally; that is to say that if $i<_P j$ then $i<j$ in the order of natural numbers. Note that this does not pose a restriction on the results, it only makes them easier to state.   Denote by $\H(P)$ the graph of the Hasse diagram of $P$. The  directed transitive closure of a graph $H$ is denoted by $\overline{H}$, and it is the graph on vertex set $V(G)$ with edges $(i,j) \in \overline{H}$ if there is an increasing path from $i$ to $j$ in $H$.

\subsection{$\G$ in terms of alternating posets.} This subsection explains how to reduce the computation of $\G$ to the computation of $\G$ for posets $P$ whose Hasse diagram is an alternating graph. A graph $G$ on the vertex set $[n]$ is called \textbf{alternating}, if there are no edges $(i,j)$ and $(j,k)$ in it with $i<j<k$. We call a poset on $[n]$ an \textbf{alternating poset} if its Hasse diagram is an alternating graph.

\bp \label{p:lr} For any naturally labeled poset $P$ on $[n]$ we can write \be \label{eq:lr} \G=\sum_{L,R} \Psi_{P_{L,R}}(\x), \ee
where the summation runs over all $L, R$ such that $L\sqcup R=[n]$,   and   $G_{L,R}=([n], \{(i,j) \in E(G) \mid i \in L, j \in R, i<j \})$ is a connected graph, where $G=\overline{\H(P)}$. Furthermore, $\H(P_{L,R})=G_{L,R}$ for a naturally labeled poset $P_{L,R}$.
\ep

\proof Recall that $\G=s(\K; \x)$. If $\K=\bigcup_{i=1}^l K_i$ for interior disjoint cones $K_i$, $i \in [l]$, then we have that $s(\K;\x)=\sum_{i=1}^l s(K_i; \x)$. If $K_i=K_{P_i}^{\rm root}$ for some posets $P_i$, $i \in [l]$, then we get that $ \G=\sum_{i=1}^l \Psi_{P_{i}}(\x)$. Therefore, to prove Equation \eqref{eq:lr}, it suffices to show that $\K=\bigcup_{L,R} K_{P_{L,R}}^{\rm root}$, where the union runs over all $L, R$ such that $L\sqcup R=[n]$,   $G_{L,R}$ is a connected 
graph ($G=\overline{\H(P)}$) and $\H(P_{L,R})=G_{L,R}$ for a naturally labeled poset $P_{L,R}$.

Since $K_P^{\rm root}=\R_+\{e_i-e_j \mid i<_P j\}$, if  $\tilde{Q}_G=\bigcup  \tilde{Q}_{G_{L,R}}$ ($\tilde{Q}_{G_{L,R}}$'s are interior disjoint), where the union runs over all $L, R$ such that $L\sqcup R=[n]$, and  $G_{L,R}$ is a connected graph, then we also obtain that $\K=\bigcup_{L,R} K_{P_{L,R}}^{\rm root}$ for interior disjoint cones $K_{P_{L,R}}^{\rm root}$. The equation $\tilde{Q}_G=\bigcup  \tilde{Q}_{G_{L,R}}$  follows from \cite[Proposition 13.3]{p1} together with the observation that $G=\overline{G}$ for our choice of $G$. 
\qed

\medskip

We note that the cones $K_{P_{L,R}}^{\rm root}$ are generally not simplicial. One way to compute $\Psi_{P_{L,R}}(\x)$ would be to triangulate $K_{P_{L,R}}^{\rm root}$  into simplicial cones with rays of the form $e_i-e_j$, since for such a cone the following simple lemma gives the value of Greene's function.

\bl \label{uni} \cite{reiner} The cone $\K$ is simplicial if and only if the Hasse diagram of $P$ contains no cycles. In this case it is also unimodular and $$\G=\frac{1}{\prod_{i\lessdot_P j }(x_i-x_j)}.$$
\el

We remark that a different proof of Lemma \ref{uni} from that given in \cite{reiner} follows immediately using the subdivision algebra of root polytopes defined in \cite{root1}.

\medskip
 \subsection{Calculating $\G$ for an alternating poset $P$.}  In light of Proposition \ref{p:lr}, if we can calculate $\G$ for an alternating poset $P$, then we can in turn calculate $\G$ for any poset $P$. In this section we accomplish the former, building on the results of Li and Postnikov \cite{lp}.  The next paragraph follows the exposition of \cite{lp}. 
 
 Given an alternating graph $G$ on the vertex set $[n]$, 
 pick a
linear order $\mathcal{O}$ on the edges of $G.$
Let $T$ be a spanning tree  of $G$, and let $e$ be an edge that does not
belong to $T$.  Let $C$ be the unique cycle contained in the graph $([n], E(T) \cup
\{e\})$.  Let $e^*$ be the maximal edge in the cycle $C$ in the linear
ordering $\mathcal{O}$ of the edges.
We say that an edge $e$ is \textbf{externally semi-active} if  either $e=e^*$ or there is
an odd number of edges in $C$ between $e$ and $e^*$.  (Since $G$ is alternating,   all
cycles in $G$  have an even length.)
Let ${\rm ext_G^{\mathcal{O}}}(T)$ be the number of externally semi-active edges of $G$  with respect
to a spanning tree $T$.

\bt \label{thm:lp} \cite{lp} Given an alternating graph $G$ and a linear ordering $\mathcal{O}$ of its edges,  let  $\T_G^{\mathcal{O}}$ be the set of spanning trees $T$  with  ${\rm ext_G^{\mathcal{O}}}(T)=0$. Then \be \tilde{Q}_G=\bigcup_{T \in \T_G^{\mathcal{O}}} \tilde{\Delta}_T,\ee where the simplices $\tilde{\Delta}_T$ are interior disjoint.
\et

 \bc \label{genp} For any naturally labeled poset $P$ on $[n]$ we can write \be \label{eq:lr} \G=\sum_{L,R} \hspace{.1in}\sum_{T \in \T_{G_{L,R}}^{\mathcal{O}_{L,R}}} \frac{1}{\prod_{(i,j)\in E(T), i<j} (x_i-x_j)}, \ee
where the summation runs over all $L, R$ such that $L\sqcup R=[n]$,   and   $G_{L,R}=([n], \{(i,j) \in E(G) \mid i \in L, j \in R, i<j \})$ is a connected graph, where $G=\overline{\H(P)}$. Furthermore, ${\mathcal{O}_{L,R}}$ is an arbitrary linear order of the edges of $G_{L,R}$.
 \ec
 
 \proof The proof follows from Proposition \ref{p:lr}. Lemma \ref{uni} and Theorem \ref{thm:lp}. \qed
 
 \bigskip
 

\subsection{An alternative proof of the BFLR Theorem.} 
Let $P_D$ be the poset of a skew diagram $D$ as in the BFLR Theorem. Let $G_D$ be the graph $\H(P_D)$ drawn on a line with vertices from left to right: $x_r, \ldots, x_1, y_1, \ldots, y_c$ and with edges as arcs above this line. Note that the condition that $G_D$ comes from $P_D$ can be translated into the conditions that $G_D$ is bipartite on parts $\{x_1, \ldots, x_r\}$ and $\{y_1, \ldots, y_c\}$ and for each $i \in [r]$, $x_i$ is connected to $y_j$, for $j \in [a_i, b_i]$, $ i \in [r]$ where $a_1\leq \cdots \leq a_r$ and $b_1\leq \cdots \leq b_r$ and $[1,c]=\cup_{i=1}^r [a_i, b_i]$. 

Given a drawing of a graph $G$ so that its vertices $v_1, \ldots, v_n$ are arranged in this order on a horizontal line, and its edges are drawn above this line, we say that $G$ is \textbf{noncrossing} if it has no edges $(v_i, v_k)$ and $(v_j, v_l)$ with $i<j<k<l$. A vertex $v_i$ of $G$ is said to be \textbf{nonalternating} if it has both an incoming and an outgoing edge; it is called \textbf{alternating} otherwise. The graph $G$ is alternating if all its vertices are alternating.

\bl \label{gd-tri} The root polytope $\tilde{Q}_{G_D}=\bigcup_{T}\tilde{\Delta}_{T}$, where the union runs over all noncrossing alternating trees of $G_D$ and the simplices $\tilde{\Delta}_{T}$ are interior disjoint. 
\el

Since noncrossing depends on the drawing of the graph it is essential that we remember that we drew $G_D$ with vertices from left to right: $x_r, \ldots, x_1, y_1, \ldots, y_c$.  To prove Lemma \ref{gd-tri} we use the following criterion due to Postnikov \cite{p1}. 

\bl \label{crit} cf. \cite[Lemma 12.6]{p1} For two   trees $T$ and $T'$ on the vertex set $[n]$, the intersection $\tilde{\Delta}_{T} \cap \tilde{\Delta}_{T'}$ is a common face of the simplices $\tilde{\Delta}_{T}$ and $\tilde{\Delta}_{T'}$ if and only if the directed graph $$U(T, T')=([n], \{(i,j) \mid (i,j) \in E(T), i<j\}\cup \{(j,i) \mid (i,j) \in E(T'), i<j\}),$$ has no directed cycles of length at least $3$.
\el

We note that  \cite[Lemma 12.6]{p1} is stated less generally then Lemma \ref{crit}, however, Postnikov's proof of it can be adapted to prove the above statement. 
\medskip

\noindent {\it Proof of Lemma \ref{gd-tri}.} One can check that the noncrossing alternating spanning trees of $G_D$ satisfy the conditions of Lemma \ref{crit}. Furthermore, one can also check that no other alternating spanning tree of $G_D$ satisfies Lemma \ref{crit} with every single noncrossing alternating spanning tree of $G_D$. Moreover, since  $\tilde{\Delta}_{T}$ is a top dimensional simplex in some triangulation of  $\tilde{Q}_{G_D}$ if and only if $T$ is an alternating tree (see \cite[Lemma 13.2]{p1}), then we are done. \qed

\bl \label{bij}The noncrossing alternating spanning trees of $G_D$ are in bijection with the lattice paths   $\pi$ from $(1,1)$ to $(r,c)$ inside $D$ that take steps either one unit south or west.
\el 
 
\begin{figure}
\begin{center}
\includegraphics[scale=.85]{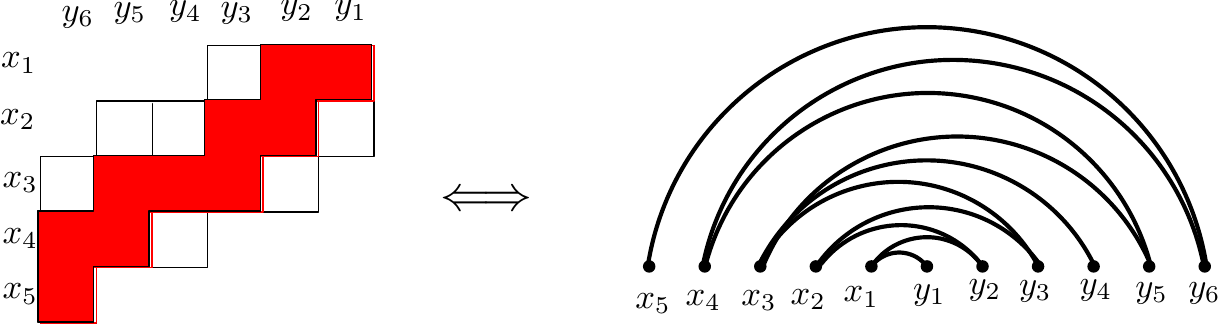}
 \caption{The correspondence between noncrossing alternating spanning trees of $G_D$ and  lattice paths    from $(1,1)$ to $(r,c)$ inside $D$ that take steps either one unit south or west. }
 \label{fig:bij}
 \end{center}
\end{figure}

 \proof The bijection is given by the map that takes a noncrossing alternating spanning tree $T=(\{x_r, \ldots, x_1, y_1, \ldots, y_c\}, \{(x_i, y_j) \mid (i,j)\in S(T)\})$ of $G_D$ to the path $\pi=S(T)$.  See Figure \ref{fig:bij}. \qed

\medskip

Given a graph $G$ on the vertex set $[n]$ such that if $(i,j)\in E(G)$ then the only increasing path from $i$ to $j$ in $G$ is the edges $(i,j)$ itself, we can define the naturally labeled poset $P_G$ to be one on the set $[n]$ with  Hasse diagram  given by (the edges of) $G$. 
 
 \bc (BFLR Theorem) {For any skew diagram $D$, }

\be \Psi_{P_D}(\x)=\sum_{\pi}\frac{1}{\prod_{(i,j)\in \pi}(x_i-y_j)},\ee
{where the sum runs over all lattice paths $\pi$ from $(1,1)$ to $(r,c)$ inside $D$ that take steps either one unit south or west.}
 \ec
 
 \proof By Lemma \ref{gd-tri} we have that the cone $K_{P_D}^{\rm root}$ is triangulated into simplicial cones $K_{P_T}^{\rm root}$, where the $T$'s run over all noncrossing alternating spanning trees of $G_D$. By Lemma \ref{bij} the latter trees are in bijection with  lattice paths   $\pi$ from $(1,1)$ to $(r,c)$ inside $D$ that take steps either one unit south or west, and thus by Lemma \ref{uni} we obtain the corollary.
 \qed

\section{Lifting Greene's theorem to the subdivision algebra}
\label{sec:lift}

The objective of this section is to generalize  Greene's theorem to a relation in the subdivision algebra of root polytopes. Subdivision algebras of root polytopes were introduced and studied in \cite{root1}, where they were used for triangulating root polytopes. Subdivision algebras were also utilized for subword complexes and flow polytopes in \cite{pipe1, pipe2, mm, prod, h-poly1, h-poly2}. We will see in this section that both Greene's theorem and an analogous one for the integer point transform of the root cone is a special case of a relation in the subdivision algebra.   

 We start by explaining how to use subdivision algebras to subdivide   root cones   $K_{P}^{\rm root}$. Since Greene's function of a poset $P$ is a valuation on a root cone $K_{P}^{\rm root}$ and we know its expression for unimodular root cones, if we triangulate $K_{P}^{\rm root}$ into unimodular root cones, then we obtain a way to calculate  Greene's function of $P$. 
  
\subsection{Root cones $\C(G)$ and their subdivisions.} We establish another notation for root cones here for ease of notation. 
For an arbitrary  loopless  graph $G$,  define the \textbf{root cone} \be \C(G):={\rm span_{\R_+}}(e_i-e_j \mid (i,j)\in E(G), i<j).\ee   In order for $\C(G)$ and $\C(H)$ to be distinct for distinct graphs $G$ and $H$, we will mostly consider \textbf{good graphs} $G$, which do not contain an edge $(i,j)$, $i<j$, if there is an increasing path other than the edge $(i,j)$ in $G$. (In particular, good graphs do not contain multiple edges.)   Given a graph $H$ let $g(H)$ be the unique good graph on the vertex set $V(H)$ such that $\C(H)=\C(g(H))$.  An important property of root cones is given in the Cone Reduction Lemma below, which can be expressed through  reduction rules on graphs as we now explain.

    The {\bf reduction rule for graphs:} Given   a  graph $G_0$ on the vertex set $[n]$ and   $(i, j), (j, k) \in E(G_0)$ for some $i<j<k$, let   $G_1, G_2, G_3$ be graphs on the vertex set $[n]$ with edge sets
  \begin{eqnarray} \label{graphs}
E(G_1)&=&E(G_0)\backslash \{(j, k)\} \cup \{(i, k)\}, \nonumber \\
E(G_2)&=&E(G_0)\backslash \{(i, j)\} \cup \{(i, k)\},  \nonumber \\
E(G_3)&=&E(G_0)\backslash \{(i, j), (j,k)\} \cup \{(i, k)\}.
\end{eqnarray}

    We say that $G_0$ \textbf{reduces} to $G_1, G_2$ and $G_3$ under the reduction rules defined by equations (\ref{graphs}).

    \begin{lemma} cf. \cite{root1} \label{reduction_lemma} \textbf{(Cone Reduction Lemma)} 
Given   a loopless  good graph $G_0$   let  $(i, j), (j, k) \in E(G_0)$ for some $i<j<k$ and $G_1, G_2$ as described by equations (\ref{graphs}).   Then  
\be \mathcal{C}(G_0)=\mathcal{C}(G_1) \cup \mathcal{C}(G_2), \ee and 
\be \mathcal{C}(G_3)=\mathcal{C}(G_1) \cap \mathcal{C}(G_2),\ee

\noindent where the cones $\mathcal{C}(G_0), \mathcal{C}(G_1), \mathcal{C}(G_2)$ are  of the same dimension and  $\mathcal{C}(G_3)$ is a facet of both   $\mathcal{C}(G_1)$ and $\mathcal{C}(G_2)$.
\end{lemma}

\proof \cite{root1} contains the proof of the above lemma in case $G$ is acyclic. A careful reading of the proof shows that the lemma holds in the case of loopless good graphs also. \qed

\medskip


\subsection{The subdivision algebra, Greene's Theorem and the integer point transform of a root cone.} In this subsection we explain the subdivision algebra and show how it yields a slick proof for Greene's theorem and its generalization. 

Observe that a graph $G$ can be encoded by the monomial $m[G]=\prod_{(i,j)\in E(G), i<j}x_{ij}$ and the reduction rule going from $G_0$ to $G_1, G_2$ and $G_3$  can be encoded  by the equation $x_{ij}x_{jk}=x_{ik}(x_{ij}+x_{jk}+\beta)$. We define the {\bf subdivision algebra} $\S$ of root polytopes as the commutative algebra generated by the variables $x_{ij}$, $1\leq i<j\leq n$, subject to the relations $x_{ij}x_{jk}=x_{ik}(x_{ij}+x_{jk}+\beta)$, for $1\leq i<j<k\leq n$. 

Let us explain the connection of the subdivision algebra to Greene's function. If we set $\beta=0$,  then the relation $x_{ij}x_{jk}=x_{ik}(x_{ij}+x_{jk})$ of $\S$ is satisfied by  $x_{ij}:=\frac{1}{x_i-x_j}$, which are the kind of terms appearing in Greene's function. If instead, we set $\beta=-1$, then the relation $x_{ij}x_{jk}=x_{ik}(x_{ij}+x_{jk}-1)$ of $\S$ is satisfied by  $x_{ij}:=\frac{1}{1-\frac{x_i}{x_j}}$. The latter will play a part in calculating the \textbf{integer point transform $\sigma_{K_{P}^{\rm root}}({\bf x})$ of the root cone} $K_{P}^{\rm root}\subset \mathbb{Z}^d$ defined as  \be \sigma_{K_{P}^{\rm root}}({\bf x}):=\sum_{{\bf m} \in  K_{P}^{\rm root}\cap \mathbb{Z}^d} {\bf x}^{\bf m}.\ee  The funcion $\sigma_{K_{P}^{\rm root}}({\bf x})$ can be seen as a finer invariant of the cone then  $\G$, as explained in \cite[Section 2.4]{reiner}. We note that in \cite{reiner} the integer point transform $\sigma_{K_{P}^{\rm root}}({\bf x})$ is denoted as $H(K_{P}^{\rm root}; {\bf X})$ and is referred to as the Hilbert series of the affine semigroup ring of the root cone. We chose to follow the more geometric name and notation of \cite[Section 3.2]{BR}.

We are now ready to prove the following generalization of Greene's theorem  via the subdivision algebra, which first  appeared in \cite{reiner}:  

\bt \label{thm:g} \cite[Corollary 8.10]{reiner} For any (connected) strongly planar poset $P$ on $[n]$ we have \be \label{g1}\sigma_{K_{P}^{\rm root}}({\bf x})=\frac{\prod_{\rho \in b(P)}(1-\frac{x_{min(\rho)}}{x_{max(\rho)}})}{\prod_{i\lessdot_{P}j}(1-\frac{x_i}{x_j})}\ee and \be \label{g2} \G=\frac{\prod_{\rho \in b(P)}(x_{min(\rho)}-x_{max(\rho)})}{\prod_{i\lessdot_{P}j}(x_i-x_j)},\ee where $\rho$ runs through all bounded regions of the Hasse diagram.
\et

\proof Let ${K_{P}^{\rm root}}=\C(G)$ for a loopless good graph $G$. Note that a root cone $\C(H)$ is unimodular if and only if $H$ is acyclic. We will use the Cone Reduction Lemma to write $\C(G)$ as a union of unimodular cones. Note that the Cone Reduction Lemma applies to loopless good graphs, and thus if we want to repeatedly apply it to the outcome cones $\C(G_i)$, $i \in [3]$, we need to apply it to $g(G_i)$, $i \in [3]$. 

Since $P$ is a connected strongly planar poset, it follows that 
$G$ is a good graph on the vertex set $[n]$ such that for every cycle $C$ of $G$ the only alternating vertices of $C$ (considered within $C$), that is vertices that have only incoming or only outgoing edges,  are its minimal and maximal vertices.  Therefore, we can apply the Cone Reduction Lemma repeatedly in such a fashion that at the end we end up with trees $T_1, \ldots, T_k$ (with $n-1$ edges), and forests $F_{n-i}^{j}$, $2\leq i \leq n-1$, $j \in I_{n-i}$ (for some index sets $I_{n-i}$), with $n-i$ edges, where $\C(T_1), \ldots, \C(T_k)$ are  unimodular cones triangulating $\C(G)$ and the $\C(F_{n-i}^{j})$'s are their intersections.

If we inspect what edges  we had to drop in the process to make sure we always apply the Cone Reduction Lemma to good graphs and obtain the acyclic graphs described in the previous paragraph, we find the following relation in the subdivision algebra: 

\be \label{1} m[G]=\prod_{\rho \in b(P)}x_{min(\rho),max(\rho)}(\sum_{T_i} m[T_i]+\sum_{F_{n-i}^{j}}{\beta}^{i-1} m[F_{n-i}^{j}]).\ee 

Note that 
\be \label{2} \sigma_{K_{P}^{\rm root}}({\bf x})=(\sum_{T_i} m[T_i]+\sum_{F_{n-i}^{j}}{(-1)}^{i-1} m[F_{n-i}^{j}])\mid_{x_{ij}=\frac{1}{1-x_i{x_j}^{-1}}}\ee

and 

\be \label{3} \G=\sum_{T_i} m[T_i]\mid_{x_{ij}=\frac{1}{x_i-x_j}}.\ee

Equations \eqref{1}, \eqref{2} and \eqref{3} together with the observations that $x_{ij}=\frac{1}{1-x_i{x_j}^{-1}}$ satisfies $x_{ij}x_{jk}=x_{ik}(x_{ij}+x_{jk}-1)$ and $x_{ij}=\frac{1}{x_i-{x_j}}$ satisfies $x_{ij}x_{jk}=x_{ik}(x_{ij}+x_{jk})$ immediately yield equations \eqref{g1} and \eqref{g2}.

\qed

 We can see equation \eqref{1} as the main theorem of this section, so we bestow it with that title:

\bt \label{thm:gg} Let $G=\H(P)$ of a naturally labeled connected strongly planar poset $P$. Then, using the notation of the proof of Theorem \ref{thm:g}, we have that  
\be m[G]=\prod_{\rho \in b(P)}x_{min(\rho),max(\rho)}(\sum_{T_i} m[T_i]+\sum_{F_{n-i}^{j}}{\beta}^{i-1} m[F_{n-i}^{j}]) \nonumber \ee holds in the subdivision algebra. 
\et

Both statements of Theorem \ref{thm:g} are   special cases of Theorem \ref{thm:gg} as shown in the proof of Theorem \ref{thm:g}.

\section{Generalizing Greene's Theorem beyond strongly planar posets}
\label{sec:gen}

In this section we will examine a special family of posets  for which Greene's function factors linearly. These posets  were first identified by Boussicault, F\'eray, Lascoux and Reiner  in \cite{reiner} who proved the aforementioned  result by studying the  affine semigroup ring of the root cone. We will give a short alternative proof   via root polytopes.  
 
The next  paragraph contains  definitions  following the exposition of \cite{reiner}. 

In a finite poset $P$, say that a triple of elements $(a, b, c)$ forms a
notch of {V}-shape (dually, a notch of $\wedge$-shape) if $a \lessdot_P b, c$ (dually, $ b, c\lessdot_P a$), and
in addition, $b, c$ lie in different connected components of the poset $P \backslash P_{²a}$ (dually,
$P \backslash P_{³a}$).
When $(a, b, c)$ forms a notch of either shape in a poset $P$, say that the quotient
poset $\overline{P} := P/ \{b = c\}$, having one fewer element and one fewer Hasse diagam edge,
is obtained from $P$ by closing the notch, and that $P$ is obtained from $\overline{P}$ by opening
a notch. 

\bt \label{notch}Let $P$ be a connected poset in which $(a, b, c)$ forms a notch, and let $\overline{P} := P/ \{b = c\}$. We assume without loss of generality that $P$ and $\overline{P}$ are naturally labeled.  Then the root polytope $\tilde{Q}_{\H(P)}$ has a triangulation with top dimensional simplices $\tilde{\Delta}_{T_1}, \ldots, \tilde{\Delta}_{T_k}$ and $\tilde{Q}_{\H(\overline{P})}$ has a triangulation with top dimensional simplices $\tilde{\Delta}_{T'_1}, \ldots, \tilde{\Delta}_{T'_k}$, where $(a,b)\in T'_i$, $(a,b), (a,c) \in T_i$, $i \in [k]$, and moreover $T_i\mid_{b=c}=T'_i$ (we ignore multiple edges).
\et

\proof The criterion of   Lemma \ref{crit} is sufficient to establish the above theorem,  since we also have that $\tilde{Q}_{\H(\overline{P})}$ has a triangulation with top dimensional simplices $\tilde{\Delta}_{T'_1}, \ldots, \tilde{\Delta}_{T'_k}$, where $(a,b)\in T'_i$, as $e_a-e_b$ is a vertex of  $\tilde{Q}_{\H(\overline{P})}$.
\qed

\medskip

When we calculate $\sigma_{K_{\overline{P}^{\rm root}}}({\bf x})$ and $\Psi_{\overline{P}}({\bf x})$   using triangulations of the root cones as implied by Theorem \ref{notch}, we immediately get:

\bc \label{86} \cite[Theorem 8.6]{reiner} When $\overline{P}$ is obtained from $P$ by closing a {V}-shaped notch $(a,b,c),$ then 
\be \sigma_{K_{\overline{P}^{\rm root}}}({\bf x})=(1-x_a{x_b}^{-1})\sigma_{K_{{P}^{\rm root}}}({\bf x})\mid_{x_b=x_c},\ee and 

\be \Psi_{\overline{P}}({\bf x})=(x_a-x_b)\G\mid_{x_b=x_c}.\ee

\ec

A consequence of Theorem \ref{notch} is the following generalization of Greene's theorem pertaining to  posets $P$ to which we can repeatedly apply the opening notch operation and obtain a poset whose Hasse diagrams  has only cycles as biconnected components.  Such posets $P$ we call \textbf{admissible}. We now recall the definition of biconnected components following \cite{reiner}. Given a graph $G=(V, E)$ we say that two edges of it are cycle-equivalent if there is a cycle which contains both edges. Let $E_i$ be the equivalence classes of this relation. Let $V_i$ be the set of vertices which are at least the endpoint of one edge in $E_i$. Then the {biconnected components} of $G$ are the graphs $G_i=(V_i, E_i)$. 

\bt Let $P$ be  an admissible  planar poset.  Then,  we have \be \label{u1}\sigma_{K_{P}^{\rm root}}({\bf x})=\frac{\prod_{\rho \in b(P)}(1-\prod_{i \in min(\rho)} {x_{i}}\prod_{j \in max(\rho)}{x_{j}^{-1}})}{\prod_{i\lessdot_{P}j}(1-{x_i}{x_j}^{-1})}\ee and \be \label{u2} \G=\frac{\prod_{\rho \in b(P)}(\sum_{i \in min(\rho)}x_{min(i)}-\sum_{j \in max(\rho)}x_{j})}{\prod_{i\lessdot_{P}j}(x_i-x_j)},\ee where $\rho$ runs through all bounded regions of the Hasse diagram of $P$.
\et

\proof This theorem can be deduced from Corollary \ref{86} together with Corollaries 8.2 and 8.3 appearing in \cite{reiner}. We note that the latter Corollaries also have simple proofs using the root polytope considerations of this paper, and we leave such alternative proofs as an exercise for the interested reader.
\qed

\section*{Acknowledgements}
I am grateful to Vic Reiner for bringing  Greene's function to my attention as well as for many  informative and valuable  exchanges about this work. I am also grateful to Alex Postnikov for  sharing his knowledge generously.

\bibliography{biblio-kir}
\bibliographystyle{plain}

\end{document}